In what follows we will consider periodic SMP matrices.
These were introduced and partially studied in \cite{smp}.
Periodic SMP matrices can be naturally extended to the negative half axis. Similarly
to the case of Jacobi matrices, such an extension appears to be very useful.

Following the result of the 
previous theorem we make the following definition:

\begin{defsmp}
 \label{def:smp}
We say that a two-sided real self-adjoint pentadiagonal matrix $A$ is \emph{SMP-structured} if its entries are uniformly
bounded and, additionally, if all even entries
on the most outter diagonals vanish and all the odd ones are positive. This means that $A$ can be written in the 
following way
\begin{equation}
\label{eqn:defsmp}
  A = S^2 \br + S \bp + \bq + \bp S^{-1} + \br S^{-2} ,
\end{equation} 
where  $\bp = \diag\{p_n\}$, $\bq = \diag\{q_n\}$, $\br = \diag\{r_{n}\}$ are real diagonal matrices 
with $r_{2n+1} > 0$, $r_{2n} = 0$ and where the sequences
$\{p_n\}$, $\{q_n\}$ as well as  $\{r_{n}\}$ are uniformly bounded. Furthermore, we demand that $\{r_{2n+1}\}$
is uniformly bounded from zero, i.e. there exists an $\eta>0$ such that $r_{2n+1}\geq\eta$ for all $n\in\mathbb Z$.

We now call $A$ an \emph{SMP matrix}, or $A\in\smp$, if it is invertible and both $A$ and $A^\tau\deq-SA^{-1}S^{-1}$
are SMP-structured.
\end{defsmp}

Similar to Jacobi matrices we can represent SMP matrices as a two-dimensional pertubation of a block-diagonal matrix
\begin{equation*}
 A=\left(\begin{array}{c|c}
    A_{-} &\mathbf0\\
    \hline
    \mathbf0& A_+
   \end{array} \right)
+ a_0\left(\tilde{e}_0\langle\cdot,e_{-1}\rangle+e_{-1}\langle\cdot,\tilde{e}_0\rangle\right),
\end{equation*}
where $a_0=\sqrt{p_0^2+r_1^2}$ and $\tilde{e}_0=\frac{1}{a_0}(p_0e_0+r_1e_1)$. In this notation the matrix $A_{-}$ corresponds to a 
one-sided SMP matrix which was defined in Section~\ref{sec:onesided}.
\begin{notationainv}
\label{notation:ainv}
 For simplicity, in what follows, we write the inverse of an SMP matrix $A$ as
\begin{equation*}
 A^{-1}= S^2 \mathbf P + S \mathbf \Pi + \mathbf \Sigma + \mathbf \Pi S^{-1} + \mathbf P S^{-2},
\end{equation*}
 where $\mathbf P = \diag\{\rho_{n}\}$, $\mathbf \Sigma = \diag\{\sigma_n\}$ and $\mathbf \Pi = \diag\{\pi_n\}$. Hence, considering the
above definition, we have that $\rho_{2n+1}=0$ and $\rho_{2n}<0$.
\end{notationainv}

\begin{thmfreeparameters}
 \label{thm:freeparameters}
Let $r_{2n}=0$, $r_{2n+1}>0$ and let the sequences $\{\log r_{2n+1}\}$, $\{p_n\}$, $\{q_{2n+1}\}$ be uniformly bounded. 
Then,
$A$ is an SMP matrix if and only if  there exists an $\varepsilon>0$ such that
\begin{equation}
\label{eqn:dn}
 q_{2n-1} p_{2n-2} p_{2n+1} - p_{2n-2} p_{2n} r_{2n+1} - p_{2n-1} p_{2n+1} r_{2n-1}\leq -\varepsilon
\end{equation}
and the coefficient sequences satisfy
\begin{equation*}
 q_{2n}=\frac{p_{2n} p_{2n+1}}{r_{2n+1}}.
\end{equation*}

Moreover, in this case, $A^{-1}$ exists and its elements  are given by
\begin{IEEEeqnarray*}{rCl}
 \rho_{2n+1}&=&0,\\
 \rho_{2n} &=& \frac{r_{2n-1} r_{2n+1}}{q_{2n-1} p_{2n-2} p_{2n+1} - p_{2n-2} p_{2n} r_{2n+1} - p_{2n-1} p_{2n+1} r_{2n-1}},\\
\pi_{2n+1} &=& -\frac{p_{2n+3} \rho_{2n+2}}{r_{2n+3}},\\
\pi_{2n} &=& -\frac{p_{2n-2} \rho_{2n}}{r_{2n-1}}, \\
\sigma_{2n+1} &=& \frac{p_{2n} p_{2n+3} \rho_{2n+2}}{r_{2n+1} r_{2n+3}} $~and$\\
\sigma_{2n} &=& -
\frac{p_{2n}(\pi_{2n+1}r_{2n+1}+q_{2n-1}\pi_{2n}+p_{2n-1}\rho_{2n})}{p_{2n}^2+p_{2n+1}^2}\\
&&-\frac{p_{2n+1}(\pi_{2n}r_{2n+1}+q_{2n+1}\pi_{2n+1}+p_{2n+2}\rho_{2n+2})}{p_{2n}^2+p_{2n+1}^2}.
\end{IEEEeqnarray*}
\end{thmfreeparameters}

\begin{proof}
 First, assume that $A\in\smp$. From $A A^{-1}=\id$ we can derive the following identities
\begin{IEEEeqnarray*}{c}
\phantomsection
\label{eqn:atimesainv}
 S^4 \left[(S^{-2} \br S^2) \mathbf P \right] =\mathbf0,\quad \left[\br(S^{-2}\mathbf P S^2)\right]S^{-4}=\mathbf0, 
\\
S^3 \left[( S^{-1} \br S) \mathbf \Pi + (S^{-2} \bp S^2) \mathbf P \right] = \mathbf0, \quad
\left[\bp(S^{-1}\mathbf P S)+\br(S^{-2}\mathbf\Pi S^2)\right]S^{-3}=\mathbf0,
\\
S^2 \left[ \br \mathbf{\Sigma} + (S^{-1} \bp S) \mathbf \Pi +(S^{-2} \bq  S^2) \mathbf P \right] = \mathbf0,
\\
\left[\bq\mathbf P+\bp(S^{-1}\mathbf \Pi S)+\br(S^{-2}\mathbf\Sigma S)\right]S^{-2}=\mathbf0,
\IEEEyesnumber\\
S \left[ S \br  \mathbf \Pi S^{-1} + \bp \mathbf \Sigma +
(S^{-1} \bq S) \mathbf \Pi + (S^{-1} \bp S) \mathbf P \right] =\mathbf0,
\\
\left[S\bp \mathbf P S^{-1}+\bq\mathbf\Pi+\bp(S^{-1}\mathbf \Sigma S)+\br(S^{-1}\mathbf\Pi S)\right]S^{-1}
=\mathbf0,
\\
S^2 \br  \mathbf P S^{-2} + S \bp \mathbf \Pi S^{-1} + \bq \mathbf \Sigma +
\bp \mathbf \Pi + \br \mathbf P = \id,
\end{IEEEeqnarray*}
where $\mathbf0$ denotes the zero operator.

From the first equation we immediately see that $\rho_{2n+1}=0$. Considering the fact that, by definition, $r_{2n+1}$ is uniformly
greater than zero, we obtain
\begin{equation*}
 \pi_{2n+2}=-\frac{p_{2n}\rho_{2n+2}}{r_{2n+1}}
\end{equation*}
from the third equation. Inserting this into the equation for the second upper diagonal yields
\begin{equation*}
 q_{2n}=-\frac{p_{2n+1}\pi_{2n+2}}{\rho_{2n+2}}=\frac{p_{2n}p_{2n+1}}{r_{2n+1}},
\end{equation*}
where we used that $\rho_{2n}<0$ uniformly. We may now continue by solving the above system for the entries of $A^{-1}$ resulting
in the identities from the claim. Since the inverse of an invertible matrix is unique, these are the coefficient sequences of 
$A^{-1}$ indeed. We once again refer to the fact that both $r_{2n+1}$ and $-\rho_{2n}$ are uniformly
bounded from above and from zero. Together with the expression which we obtained for $\rho_{2n}$ this implies (\ref{eqn:dn}).

For the proof of the \emph{if}-part, we first of all verify that all the sequences involved are uniformly bounded. For $q_{2n}$
this, again, follows from the fact that so is $\log r_{2n+1}$. Considering (\ref{eqn:dn}) as well yields that this claim holds
true for $\log|\rho_{2n}|$. Hence, it also does for $\pi_{2n}$, $\pi_{2n+1}$ and $\sigma_{2n+1}$. 
In addition to that, we see from (\ref{eqn:dn}) that $p_{2n}^2+p_{2n+1}^2$, too, is uniformly greater than zero. Indeed,
suppose that w.l.o.g. $p_{2n}$ and $p_{2n+1}$ converge to zero. Then the left handside of
(\ref{eqn:dn}) tends to zero as $n\to\infty$, which clearly leads to a contradiction.
Consequently, $\sigma_{2n}$ is uniformly bounded too.

Secondly, note that the identities from the claim fulfill all the equations from (\ref{eqn:atimesainv}). Thus, $A$ is invertible
and, additionally, the coefficient sequences of its inverse satisfy
the conditions from Definition~\ref{def:smp} and hence $A\in\smp$.
\end{proof}